\documentclass[a4paper,12pt]{article}
\usepackage{amsfonts}
\usepackage{amsthm}
\usepackage{amsmath}
\usepackage{amsfonts}
\usepackage{latexsym}
\usepackage{amssymb}
\usepackage[latin1]{inputenc}

\newcommand{\U}{\mathcal{U}}
\newcommand{\um}{\U_{\cM}}

\def\a{\mathcal A}

\def\P{\mathcal P}

\newcommand{\CC}{\mathbb{C}}
\newcommand{\RR}{\mathbb{R}}
\newcommand{\NN}{\mathbb{N}}

\newtheorem{fed}{Definition}[section]
\newtheorem{teo}[fed]{Theorem}
\newtheorem*{teo*}{Theorem}
\newtheorem{lem}[fed]{Lemma}
\newtheorem{cor}[fed]{Corollary}
\newtheorem{pro}[fed]{Proposition}
\theoremstyle{definition}
\newtheorem{rem}[fed]{Remark}

\def\bmat{\left[\begin{array}}
\def\emat{\end{array}\right]}
\def\ben{\begin{enumerate}}
\def\een{\end{enumerate}}

\def\co{\mathrm{conv}}

\def\M{\mathcal{M}}
\def\um{\mathcal{U}_\M}
\newcommand{\uni}[1]{\um( #1 )}
\newcommand{\cluni}[1]{\overline{\um( #1 )}}

\def\cA{\mathcal{A}}

\def\cP{\mathcal{P}}

\newcommand{\hil}{\mathcal{H}}
\newcommand{\op}{B(\mathcal{H})}

\newcommand{\peso}[1]{\quad \text{ #1 } \quad }

\begin{document}
\title{Refinements of spectral resolutions and
modelling of operators in II$_1$ factors}
\author{Pedro G. Massey\footnote{Supported in part by Consejo Nacional de Investigaciones Científicas y Técnicas of Argentina}}
\date{}
\maketitle

\begin{center}
{\it A Marina, Candelaria y Agustín, con amor} \end{center}

\begin{abstract}

We study refinements between spectral resolutions in an arbitrary
II$_1$ factor $\M$ and obtain diffuse (maximal) refinements of
spectral resolutions. We construct models of operators with
respect to diffuse spectral resolutions. As an application we
obtain new characterizations of sub-majorization and spectral
preorder between positive operators in $\M$ and new versions of
some known inequalities involving these preorders.
\end{abstract}

\noindent {\small Keywords:  II$_1$ factor, bounded right spectral
resolution, spectral preorder, sub-majorization.\\ AMS subject
classification (2000) Primary 46L51; Secondary 47A63.}

\section{Introduction}

The study of the norm closure of unitary orbits of self-adjoint
operators in von Neumann algebras is a well established area of
research. Some of the early results on this subject go back to the
work of Weyl and von Neumann in the type I factor case. Kamei, in
his development of majorization between operators in II$_1$
factors, obtained an interesting characterization of the norm
closure of the unitary orbit of a positive operator in terms of
its singular values. Recently, Arveson and Kadison have described
these sets for self-adjoint operators in terms of spectral
distributions \cite{arvkad} in the II$_1$ factor and  Sherman
\cite{she} has obtained interesting descriptions of several
closures of unitary orbits in von Neumann algebras under weak
restrictions (see the introduction of \cite{she} for a detailed
account on the history of these problems and recent references).
It turns out that even in the general setting of \cite{she}, the spectral data
of operators play a fundamental role in these investigations.

There are other notions closely related to unitary orbits, that
are defined in terms of spectral data, such as majorization,
sub-majorization and spectral dominance; the study of these
notions has been considered in several research works like the
papers of Kamei \cite{kam} and Hiai \cite{Hiai0,Hiai}, Hiai and
Nakamura \cite{HiaiN0,{HiaiN}} and the more recent papers of
Kadison \cite{kad1,kad2,kad} and of Arveson and Kadison
\cite{arvkad}. In this context one usually tries to describe
operators in some set associated with (the norm closure of)
$$\U_\M(b):=\{u^*bu:\ u\in \M \ \text{is a unitary operator}\}$$
where $\M$ is a semifinite von Neumann algebra with faithful
semifinite trace $\tau$ and $b\in \M$ is a self-adjoint operator.
For example, it is well known \cite{kam} that if $\M$ is a II$_1$
factor then $a\in\co(\overline{\U_\M(b)})$ if and only if $a$ is
majorized by $b$, which is a spectral relation. In this case the
spectral data of $a$ may be more complex (disordered) than
that of $b$. This makes things difficult when trying to recover
$a$ as an element of $\co(\overline{\U_\M(b)})$ whenever we know
that $a$ is majorized by $b$. In order to overcome a similar
difficulty, in \cite{argmas} we considered an ``diffuse"
refinement of the (joint) spectral measure of an ordered $n$-tuple
of mutually commuting self-adjoint elements of a II$_1$ factor
$\M$.

In this work we consider a related construction to that obtained
in \cite{argmas} that, roughly
speaking, allows us to represent every positive operator $a\in
\M^+$ as Borel functional calculus (by an increasing
left-continuous function) of a positive operator $a'\in \M$ with
maximal disordered spectral resolution (with respect to a
preorder called \emph{refinement} that we shall introduce).
Moreover, the operator $a'\in \M^+$ has the following property:
\emph{any} positive operator $b\in \M^+$ is, up to approximately
unitary equivalence, Borel functional calculus of $a'$ (by an
increasing left-continuous function). These constructions are what
we call \emph{diffuse refinements of spectral resolutions} and
\emph{modelling of operators}. We also consider some relations
between these constructions and maximal abelian subalgebras of
$\M$. The idea of considering maximal (diffuse) refinements of
spectral resolutions and of constructing some models of operators
in finite factors has already been considered in
\cite{HiaiN0,{HiaiN}} although the notion of refinement introduced
here has not. In this work we attempt a brief but systematic
treatment of these concepts.

Our results are related to Kadison's study of Schur-type
inequalities \cite{kad} and Arveson-Kadison's study of closed
unitary orbits in II$_1$ factors \cite{arvkad}. Indeed our
techniques provide alternative proofs to some of their results.
Moreover, our refinements and modelling techniques are the basis
for a version of the Schur-Horn type theorem in II$_1$
factors in  \cite{argmas2}.

 As an application of these constructions we
present characterizations of the sets $$\{c\in
\M:\ 0\leq c\leq d\in \overline{\U_\M(a)}\}$$ and $$\{c\in \M:\
0\leq c\leq d\in \overline{\co(\U_\M(a)})\}$$ in terms of spectral
data. These characterizations are then applied to some recent
spectral inequalities obtained in \cite{JD,{silva},{Doug}}.

 The paper is organized as follows. In
section 2 we recall some definitions and facts regarding spectral
relations (spectral preorder, majorization and sub-majorization).
In section 3 we present our results on refinements of bounded
right spectral resolutions in II$_1$ factors. In section
4 we consider the modelling of operators and use this construction
to study spectral dominance and sub-majorization.

\section{Preliminaries}

Let $\op$ be the algebra of bounded operators on a Hilbert space
$\hil$. In what follows, the pair $(\M,\tau)$ shall denote a
semifinite von Neumann algebra and a faithful normal semifinite
(f.n.s.) trace on $\M$. In particular, if $\M$ is a finite factor
then $\tau$ denotes the unique f.n.s. trace such that $\tau(1)=1$.
The real space of self-adjoint operators in $\M$ is denoted by
$\M_{sa}$, the cone of positive operators by $\M^+$ and the
unitary group by $\U_\M$. If $a\in \M_{sa}$ then $P^a(\Delta)$
denotes the spectral projection of $a$ corresponding to the
measurable set $\Delta\subseteq \RR$. For simplicity of notation
we shall write $P^a(\alpha,\beta)$ (instead of $
P^a((\alpha,\beta))$) for a real interval $(\alpha,\beta)\subseteq
\RR$. $\cP(\M)\subseteq \M_{sa}$ denotes the lattice of orthogonal
projections in $\M$ endowed with the strong operator topology. For
$a\in \M$, $R(a)$ denotes its range and $P_{\overline{R(a)}}\in
\cP(\M)$ the orthogonal projection onto the closure of its
range. By a decreasing function (resp. increasing) we mean a
non-increasing function (resp. non-decreasing). If $(X,\nu)$ is a
measure space then $L^\infty(\nu)^+$ denotes the cone of
$\nu$-essentially bounded nonnegative functions on $X$. The set of
nonnegative numbers is denoted by $\RR^+_0$.

\subsection{Singular values, spectral preorder and (sub) majorization}

The \emph{$\tau$-singular values} (or $\tau$-singular numbers)
\cite{fack} of $x\in \M$ are defined for each $t\in \RR^+_0$ by
\begin{equation}\label{ecua prelis1}
\mu_x(t)=\inf\{\|xe\|:\ e\in \mathcal P(\M),\,\tau(1-e)\leq
t\}.\end{equation} The function $\mu_x:\RR^+_0\rightarrow \RR^+_0$
is  decreasing and right-continuous. If $x,\,y\in \M$ then
\begin{equation}\label{ecua prelis2} |\mu_x(t)-\mu_y(t)|\leq \|x-y\|
\end{equation}
 which shows a continuous
dependence of the singular values on the operator norm. If $a\in
\M^+$, we have $$\mu_a(t)=\min\{s\in \RR^+_0:\
\tau(P^a(s,\infty))\leq t\}.$$ This last characterization of the
singular values of positive operators shows the following
property: if $a,\,b\in \M^+$ are such that
$\tau(P^a(s,\infty))=\tau(P^b(s,\infty))$ for every $s\in \RR^+_0$
then $\mu_a=\mu_b$. On the other hand, from \eqref{ecua prelis1}
we see that $\mu_a=\mu_{uau^*}$ for every unitary operator $u\in
\U_\M$. Moreover, from this last fact and the continuous
dependence \eqref{ecua prelis2} we see that $\mu_a=\mu_b$,
whenever $b\in \overline{\U_\M(a)}$, where $\overline{\U_\M(a)}$
denotes the norm closure of the unitary orbit $$\uni{a} = \{ uau^*
: u \in \um\}.$$ Kamei proved \cite{kam} a converse of this fact
when $(\M,\tau)$ is a finite factor. We summarize these remarks in
the following proposition.
\begin{pro}\label{rems}
Let $(\M,\tau)$ be a semifinite von Neumann algebra and let
$a,\,b\in \M^+$.
\begin{enumerate}
\item[1.] If $b \in \cluni{a}$, then
$\mu_a=\mu_b$.
\item[2.] (Kamei \cite{kam}) Assume further that $(\M,\tau)$ is a finite factor and
$\mu_a=\mu_b$. Then $b \in \cluni{a}$.
\end{enumerate}
\end{pro}

Next we recall the definitions of three different preorders that
we shall consider in the sequel. If $a,\,b\in \M^+$ we say that
$b$ \emph{spectrally dominates} $a$, and write $a\precsim b$, if
any of the following (equivalent) statements holds:
\begin{enumerate}
\item[a)] $\mu_a(t)\leq \mu_b(t)$, for all $t\geq 0$.
\item[b)] $\tau(P^a(t,\infty))\leq\tau(P^b(t,\infty))$, for all $t\geq 0$.
\end{enumerate}
If in addition $(\M,\tau)$ is a semifinite factor
\begin{enumerate}
\item[c)] $P^a(t,\infty)\precsim P^b(t,\infty)$ in the Murray-von Neumann's sense.
\end{enumerate}
We say that $a$ is \emph{sub-majorized} by $b$, and write
$a\prec_w b$, if $$ \int_0^s \mu_a(t)\ dt\leq \int_0^s \mu_b(t) \
dt,\ \ \ \text{for every }s\geq 0. $$  If in addition
$\tau(a)=\tau(b)$ then we say that $a$ is \emph{majorized} by $b$
and write $a\prec b$. 
 It is well known that $a\leq b\Rightarrow a\precsim
b\Rightarrow a\prec_w b$.

We shall need the following result due to Hiai and Nakamura
 \cite{HiaiN}, concerning functions in a finite
 measure space $(X,\nu)$. In this
 case, a function $g\in L^\infty(\nu)$ is considered as an
 operator in the finite von Neumann algebra $(L^\infty(\nu),\varphi)$
  and singular values are defined with
 respect to the normal faithful finite trace $\varphi$ induced by $\nu$, i.e.
 \begin{equation}\label{revisado1}
 \varphi(g):=\int_{X} g\ d\nu, \ \ g\in
 L^\infty(\nu).\end{equation}
 \begin{pro}\label{pro HiaiNak}Let $(X,\nu)$ be a probability
 space and
 let $f,\,g\in L^\infty(\nu)^+$. Then
$f\prec_w g$ if and only if there exists $h\in L^\infty(\nu)^+$
such that $f\leq h\prec g$.
\end{pro}

\begin{rem}\label{fundamental sobre funciones}
If $(\M,\tau)$ is a finite factor and $a\in\M^+$, then let $\nu$
be the regular Borel probability measure given by $\nu
(\Delta)=\tau(P^a(\Delta))$. For every $g\in L^\infty(\nu)^+$ let
$$g(a)= \int_{\sigma(a)} g\ dP^a\in \M^+$$ and note that
$\mu_{g(a)}=\mu_g$. As a consequence we get that
$\tau(g(a))=\varphi(g)$, where $\varphi$ is given by
\ref{revisado1}. Thus, if $h,\,g\in L^\infty(\nu)^+$, then
$h(a)\precsim g(a)$ (resp. $h(a)\prec g(a)$, $h(a)\prec_w g(a)$)
in $\M$ if and only if $h\precsim g$ (resp. $h\prec g$, $h\prec_w
g$) in $L^\infty(\nu)$. \qed
\end{rem}


\section{Refinements of spectral resolutions}

Let $I=[\alpha,\beta]\subseteq \RR$ be a closed interval, and
recall that $\cP(\M)$ denotes the lattice of orthogonal
projections in $\M$ endowed with the strong operator topology. If
$p\in \cP(\M)$, we say that a map $E:I\rightarrow \cP(\M)$ is a
\emph{bounded right spectral resolution of} $p$ (abbreviated
``brsr of $p$") if $E$ is decreasing and right-continuous,
$E(\beta)=0$ and $E(\alpha)=p$. If $p=1$ then this notion agrees
with the usual definition of brsr in $\M$. For example, any $a\in
\M^+$ induces a brsr of $p=P_{\overline{R(a)}}\,$, by
\begin{equation}\label{la resu de a}
E(\lambda)=P^a(\lambda,\,\infty),\ \ \lambda\in [0,\|a\|].
\end{equation} Given $E:I\rightarrow \P(\M)$
 a brsr (of $E(\alpha)$) then, we identify $E$ with the family
$\{E_\lambda\}_{\lambda\in I}$, where $E_\lambda=E(\lambda)$ for
every $\lambda\in I$. If the set $I$ is clear from the context, we
simply write $\{E_\lambda\}$.

If $E:[\alpha,\beta]\rightarrow \cP(\M)$ is a brsr, we say that
$\lambda_0\in (\alpha,\beta]$ is an atom for $\{ E_\lambda\}$, if
the resolution is not continuous at $\lambda_0$; if $p\neq 1$ then
$\alpha$ is considered as an atom. The set of atoms of
$\{E_\lambda\}$ is denoted by At$(\{E_\lambda\})$. We say that
$\{E_\lambda\}$ is a \emph{diffuse} brsr if the set
At$(\{E_\lambda\})$
 is empty. It is clear that $\{E_\lambda\}$ is
 diffuse if and only if $E(\alpha)=1$ and $E$ is a continuous
 function (recall that $\cP(\M)$ is endowed with the SOT).
We say that a positive operator $a\in \M^+$ has \emph{continuous
distribution} if the resolution induced by $a$ (see (\ref{la resu
de a})) is diffuse. Therefore, $a\in \M^+$ has continuous
distribution if and only if $P_{\overline{R(a)}}=1$ and
$P^a(\{x\})=0$ for every $x\in \RR$.

It is well known that given a brsr $\{E_\lambda\}_{\lambda\in I}$
 in $\M$ then there exists a unique spectral measure $F$ on
$I$ with values in  $\mathcal P(\M)$ such that
$E_\lambda=F((\lambda,\infty))$ for every $\lambda\in I$. If
$h:I\rightarrow \CC$ is a uniformly bounded measurable function
then we use the following notation
\begin{equation}\label{int contra resu}
\int_I h(\lambda)\ dE_\lambda:=\int_I h\ dF.\end{equation}

\begin{fed}\label{defi de refinamiento}\rm
Let $\{E_\lambda\}_{\lambda\in I}$ and $\{E'_\lambda\}_{\lambda\in
I'}$ be brsr's, where $I=[\alpha,\beta]$ and
$I'=[\alpha',\beta']$. We say that $\{E'_\lambda\}$ \emph{refines}
$\{E_\lambda\}$ if there exists $f:I\rightarrow I'$ such that \ben
\item[(a)] $f$   is increasing, right-continuous and $f(\beta)=\beta'$;
\item[(b)] $E_\lambda=E'_{f(\lambda)}$ for every $\lambda\in I$.
\een We say that $\{E'_\lambda\}$ \emph{strongly refines}
$\{E_\lambda\}$ if $f$ also satisfies
 \ben
\item [(c)]
$f(\lambda)\geq \lambda$ for every $\lambda\in I$, and
\item [(d)] $f(\lambda)-f(\mu)\geq \lambda-\mu$, for every
$\lambda>\mu\in I$.  \een
\end{fed}

If $\{E'_\lambda\}$ (strongly) refines $\{E_\lambda\}$ we also say
that $(\{E'_\lambda\},f)$ is a {(strong) refinement} of
$\{E_\lambda\}$, where $f$ is as in Definition \ref{defi de
refinamiento}. It is easy to see that refinement is a preorder
relation.

The following, which is the main result of this section, is
related with the refinement of spectral measures of separable
abelian $C^*$-subalgebras in a II$_1$ factor developed in
\cite{argmas}.

\begin{teo}\label{prepo} Let $(\M,\tau)$ be a II$_1$ factor and let $a\in
\M^+$. Then there exists $a'\in \M^+$ with continuous distribution
and such that the brsr induced by $a'$ strongly refines the brsr
induced by $a$. Further, if $a\in \cA^+$, where $\cA$ is a masa in
$\M$, then $a'$ can be selected from $\cA$.
\end{teo}

In what follows we state some lemmas and use them to prove Theorem
\ref{prepo} at the end of this section. In the rest of the paper,
the pair $(\M,\tau)$ will always denote a II$_1$ factor. Let
$I=[\alpha,\beta]$ and let $\{ E_\lambda\}_{\lambda\in I}$ be a
brsr
 of a projection $p\in\cP(\M)$. If $\lambda_0\in
(\alpha,\beta]$ is an atom for $\{ E_\lambda\}$, then
\begin{equation}\label{defi de los p}
\lim_{\lambda\rightarrow\lambda_0^-} E_\lambda
=E_{\lambda_0}+p(\lambda_0), \ \ \ p(\lambda_0)\neq 0.
\end{equation}
 In this
case $p(\lambda_0)\in \cP(\M)$ is the \emph{jump projection} of
$\{E_\lambda\}$ at $\lambda_0$. If $p\neq 1$ then $\alpha\in$
At$(\{E_\lambda\})$ and the jump projection at $\alpha$ is by
definition $p(\alpha)=1-p$.  Note that the set of atoms
At$(\{E_\lambda\})$ is countable. Indeed, if
$\lambda_0,\,\lambda_1\in$ At$(\{E_\lambda\})$ and $\lambda_0\neq
\lambda_1$, then it is easy to see that
$p(\lambda_0)\,p(\lambda_1)=0$, i.e. $p(\lambda_0)$ and
$p(\lambda_1)$ are orthogonal projections. Therefore
\begin{equation}\label{hay a lo sumo numerables}
\mathcal J(\{E_\lambda\}):=\sum_{\lambda\in \rm{At}(\{
E_\lambda\})} \tau(p(\lambda))= \tau\left(\sum_{\lambda\in
\rm{At}( \{ E_\lambda\})} p(\lambda)\right)\leq 1
\end{equation}
and this implies that At$(\{E_\lambda\})$ is countable. The real
number $\mathcal J(\{E_\lambda\})$ is called the \emph{total jump}
of the resolution.

\begin{lem}\label{decrece el salto total}
Let $\{E_\lambda\}_{\lambda\in I},\,\{E'_\lambda\}_{\lambda\in
I'}$ be brsr's in $\M$. If $\{E'_\lambda\}$ refines
$\{E_\lambda\}$ then $\mathcal J(\{E_\lambda\})\geq \mathcal
J(\{E'_\lambda\})$.
\end{lem}

\begin{proof}
Let $\lambda_0\in \text{At}(\{E'_\lambda\})$ and consider
$\mu_0=\min\{\mu\in I:\ f(\mu)\geq \lambda_0\}$ which is well
defined by (a) in Definition \ref{defi de refinamiento}. Then by
definition of $\mu_0$,  $f(\mu_0)\geq \lambda_0$ and
$f(\mu)<\lambda_0$ if $\mu<\mu_0$. So $$ \lim_{\mu\rightarrow
\mu_0^-}E_\mu-E_{\mu_{0}}=\lim_{\mu\rightarrow
\mu_0^-}E'_{f(\mu)}-E'_{f(\mu_0)}\geq \lim_{\lambda\rightarrow
\lambda_0^-}E'_{\lambda}-E'_{\lambda_0}\neq0,$$ since $\lambda_0$
is an atom of $\{E'_\lambda\}$. Therefore $\mu_0\in I$ is an atom
of the resolution $\{E_\lambda\}$ and we have
\begin{equation}\label{cua y med}
 \lim_{\mu\rightarrow \mu_0^-}\tau(E_{\mu})=
\lim_{\mu\rightarrow \mu_0^-}\tau(E'_{f(\mu)})>
\tau(E'_{\lambda_0})\geq \tau(E'_{f(\mu_0)})=\tau(E_{\mu_0}),
\end{equation} since $f(\mu)\rightarrow \lambda_1^{-}\leq \lambda_0$ when
$\mu\rightarrow \mu_0^-$ and $\lambda_0\in$ At$(\{E'_\lambda\})$.
We consider the following relation in At$(\{E'_\lambda\})$: if
$\lambda_1,\lambda_2\in \text{At}(\{E'_\lambda\})$ then
$\lambda_1\approx \lambda_2$ if and only if there exists $\mu_0\in
\text{At}(\{E_\lambda\})$ such that
\begin{equation}\label{tonta mariposa}
\tau(E'_{\lambda_1}),\,\tau(E'_{\lambda_2})\in
\left[\tau(E_{\mu_0}),\lim_{\mu\rightarrow\mu_0^-}\tau(E_{\mu})\right).\end{equation}
The inequality (\ref{cua y med}) shows that this relation is
reflexive. On the other hand it is clearly symmetric. Note that if
$\mu_1<\mu_2$ then $\lim_{\mu\rightarrow \mu_2^-}\tau(E_\mu)\leq
\tau(E_{\mu_1})$ and $$ [\tau(E_{\mu_2}), \lim_{\mu\rightarrow
\mu_2^-}\tau(E_\mu)) \cap [\tau(E_{\mu_1}), \lim_{\mu\rightarrow
\mu_1^-}\tau(E_\mu))=\emptyset.$$ So, if
$\lambda_1\approx\lambda_2$ then there exists a unique $\mu_0\in $
At$(\{E_\lambda\})$ such that (\ref{tonta mariposa}) holds, so in
particular $\approx$ is an equivalence relation.
 Therefore, for any equivalence class $Q\in\Pi=\text{At}(\{E'_\lambda\})/\approx$,
  there exists a unique atom $\mu_Q\in$ At$(\{E_\lambda\})$ such that
$$\tau(E_\lambda)\in [\tau(E_{\mu_Q}),\lim_{\mu\rightarrow
\mu_Q^-}\tau(E_{\mu}))\ \ \text{ for all }\lambda\in Q.$$ Let
$\lambda_1,\ldots,\lambda_n\in Q$ with
$\lambda_1<\ldots<\lambda_n$. Then, if $p'(\lambda_i)$ is the jump
projection of the resolution $\{E'_\lambda\}$ at $\lambda_i$ and
$p(\mu_Q)$ is the jump projection of the resolution
$\{E_\lambda\}$ at $\mu_Q$, we have
\begin{eqnarray*}\label{obs Demetrio}
\sum_{i=1}^n\tau(p'(\lambda_i))&=&\sum_{i=1}^n(
\lim_{\lambda\rightarrow \lambda_i^-}\tau(E'_\lambda)-
\tau(E'_{\lambda_i}) )\leq \lim_{\lambda\rightarrow
\lambda_1^-}\tau(E'_\lambda)-\tau(E'_{\lambda_n})\\ &\leq &
\lim_{\mu\rightarrow
\mu_Q^-}\tau(E'_{f(\mu)})-\tau(E'_{f(\mu_Q)})=\tau(p(\mu_Q))
\end{eqnarray*}
 Taking limit over $n$ if necessary, we get
$\sum_{\lambda\in Q}\tau(p'(\lambda))\leq \tau(p(\mu_Q))$.
Therefore $$\mathcal J(\{E'_\lambda\})=\sum_{Q\in
\Pi}\sum_{\lambda\in Q}\tau(p'(\lambda))\leq\sum_{Q\in
\Pi}\tau(p(\mu_Q))\leq \mathcal J(\{E_\lambda\})$$ where the
rearrangement is valid since we are considering series of positive
terms.
\end{proof}

We introduce the following notation in order to state Lemma
\ref{hay limite}.
\begin{fed}\label{defi de comp}
If $\{\alpha_k\}_{k\in \NN}\in \ell^1(\RR^+)$ we say that a
sequence $(\{E^k_\lambda\}_{\lambda\in I_k})_{k\in \NN}$ of brsr's
in $\M$ is $\{\alpha_k\}_{k\in \NN}$-\emph{compatible} if the
following conditions hold:
 \ben
\item[1.] $\exists \alpha,\,\beta\in \RR^+_0$ such that
$I_k=[\alpha,\beta+\sum_{i=1}^k\alpha_i]$ for every $k\in \NN$.
\item[2.] $(\{E^{k+1}_\lambda\},f_k)$ is a strong refinement of
$\{E^k_\lambda\}$ for every $k\in \NN$.
\item[3.] $f_k(\lambda)-\lambda\leq \alpha_k$, for
every $\lambda\in I_k$ and for every $k\in \NN$.
 \een\end{fed}

\begin{lem}\label{hay limite}Let $\{\alpha_k\}_{k\in \NN}\in \ell^1(\RR^+)$ and
$(\{E^k_\lambda\}_{\lambda\in I_k})_{k\in \NN}$ be
$\{\alpha_k\}_{k\in \NN}$-compatible. Then there exists a brsr
$\{E_\lambda\}_{\lambda\in I}$ in $\M$ such that $\{E_\lambda\}$
strongly refines $\{E^k_\lambda\}$, for every $k\in \NN$.
Moreover, if $\cA\subseteq \M$ is a masa and  $\{E^k_\lambda\}$ is in $\a$ for each $k\in \NN$, we can choose $\{E_\lambda\}$ also in $\a$.
\end{lem}

\begin{proof} For simplicity, we shall assume that $\alpha=0$.
The general case follows from this by reparametrization. Let
$I=[0,\beta+\sum_{i=1}^\infty \alpha_i]$ and for every $k\in \NN$
let $f_k:I_k\rightarrow I_{k+1}$ be as in Definition \ref{defi de
comp}. Note that, since $f_k(\lambda)\geq \lambda$ for $\lambda\in
I_k$ (condition (c) in \ref{defi de refinamiento}),
$$E^k_\lambda=E^{k+1}_{f_k(\lambda)}\leq E^{k+1}_\lambda.$$
Therefore, for each $\lambda\in I$ the sequence
$\{E^k_\lambda\}_{k\in\NN}$ is increasing, where we set
$E^k_\lambda=0$ if $\lambda\notin I_k$. Let us define
\begin{equation}\label{defi de la refinada}
E_\lambda=\bigvee_{k\in \NN}
E^k_\lambda=\lim_{k\rightarrow\infty}E^k_\lambda\in \cP(\M),\ \
\lambda\in I
\end{equation} where the limit is in the strong operator topology.
Note that, if $\a\subseteq \M$ is a masa and
$E^k_\lambda\in\P(\a)$ for every $k\in \NN$, then $E_\lambda\in
\a$. To see that $\{E_\lambda\}_{\lambda\in I}$ is a brsr note
first that $E_{\lambda_0}\geq E_{\lambda}$ if
$\lambda_0\leq\lambda$. Thus $\exists\lim_{\lambda\rightarrow
\lambda_0^+}E_\lambda \leq E_{\lambda_0}$. If
$\{\lambda_n\}\subseteq I$ is a decreasing sequence such that
$\lim_{n\rightarrow\infty}\lambda_n=\lambda_0$ then
\begin{eqnarray*} \tau(\lim_{n\rightarrow
\infty} E_{\lambda_n})&=&\lim_{n\rightarrow
\infty}\tau(E_{\lambda_n})=
\lim_{n\rightarrow\infty}\lim_{k\rightarrow\infty}
\tau(E^k_{\lambda_n})\\ &=&
\lim_{k\rightarrow\infty}\lim_{n\rightarrow\infty}
\tau(E^k_{\lambda_n})=\tau(\bigvee_{k\in\NN}E^k_{\lambda_0})=\tau(E_{\lambda_0})
\end{eqnarray*}
where the change of order of the iterated limits is valid since
the double sequence $\{\tau(E^k_{\lambda_n})\}_{n,k}$ is positive,
bounded and increasing in each variable. Therefore
$\lim_{\lambda\rightarrow \lambda_0^+}E_\lambda=E_{\lambda_0}$ and
$\{E_\lambda\}_{\lambda\in I}$ is a brsr.

Fix $k\in \NN$ and consider the sequence $\{u_n:I_{k}\rightarrow
I_{k+n}\}_{n\in \NN}$ of increasing right-continuous functions,
given inductively by $u_1=f_{k}$ and $u_n=f_{k+n-1}\circ u_{n-1}$
for $n\geq 2$. Then, it is easy to see that
\begin{enumerate}
\item[1.] $E^k_\lambda=E^{k+n}_{u_n(\lambda)}$,
\item[2.] $u_{n+1}\geq u_n$, $\|u_{n+1}-u_n\|_\infty\leq
\alpha_{n+k}$,
\item[3.] $u_n(\lambda)-u_n(\mu)\geq \lambda-\mu$ if $\lambda,\,\mu\in
I_k$ and $\lambda\geq \mu$.
\end{enumerate}
 Let $h_k:I_k\rightarrow I$ be the uniform limit of the
increasing sequence $\{u_n\}$. Then $h_k$ is increasing
right-continuous, $h_k(\lambda)\geq \lambda$ ($u_1=f_{k}$) and
$h_k(\lambda)-h_k(\mu)\geq \lambda-\mu$ if $\lambda>\mu\in I_k$.
Let $\lambda_0\in [0,\beta+\sum_{i=1}^k \alpha_i)$ and note that
$E^k_{\lambda_0}=E^{k+n}_{u_n(\lambda_0)}\geq
E^{k+n}_{h_k(\lambda_0)}$, since $u_n(\lambda)\leq h_k(\lambda)$.
Therefore
\begin{equation}\label{desigualdad} E^k_{\lambda_0}\geq
\lim_{n\rightarrow \infty}
E^{k+n}_{h_k(\lambda_0)}=E_{h_k(\lambda_0)}.\end{equation} To see
that equality holds in (\ref{desigualdad}) we consider
$$\lambda_n:=\min\{\lambda\in I_k:\ u_n(\lambda)\geq
h_k(\lambda_0)\}.$$ By definition we have  $u_n(\beta+\sum_{i=1}^k
\alpha_i)=\beta+\sum_{i=1}^{k+n} \alpha_i$ so $\lambda_n$ is well
defined. Further, $\lambda_n\geq \lambda_{n+1}\geq\lambda_0$,
since $\{u_n\}$ is an increasing sequence, and
$\lambda_n\rightarrow \lambda_0^+$. Indeed, if $\lambda>\lambda_0$
and $\lambda-\lambda_0=\epsilon$ then $h_k(\lambda)\geq
h_k(\lambda_0)+\epsilon$ and there exists $n\in \NN$ such that
$u_n(\lambda)>h_k(\lambda_0)$, which implies that $\lambda_0\leq
\lambda_n\leq \lambda$. Finally, we have $$ E_{h_k(\lambda_0)}\geq
E_{u_n(\lambda_n) }\geq E^{k+n}_{u_n(\lambda_n)}=E^k_{\lambda_n},\
\forall n\in \NN $$which implies that $E_{h_k(\lambda_0)}\geq
\lim_{n\rightarrow \infty}E^k_{\lambda_n}=E^k_{\lambda_0}$.
\end{proof}

\begin{lem}\label{primer refinamiento}
Let $\{E_\lambda\}_{\lambda\in [\alpha,\beta]}$ be a
brsr in $\M$. If $\lambda_0\in
\text{At}(\{E_\lambda\})$, then there exists a strong refinement
$(\{E'\}_{\lambda\in I' },f)$ of $\{E_\lambda\}$, where
$I'=[\alpha,\beta+\tau(p(\lambda_0))]$ such that
\begin{enumerate}
\item[1.] $\mathcal J(\{E'_\lambda\})=\mathcal
J(\{E'_\lambda\})-\tau(p(\lambda_0))$.
\item[2.] $f(\lambda)-\lambda\leq \tau(p(\lambda_0))$ for every $\lambda\in I$.
\item[3.] At$(\{E'_\lambda\})=f(\text{At}(\{E_\lambda\}\setminus
\lambda_0))$.
\end{enumerate}
Moreover, if $\a\subseteq \M$ is a masa and $\{E_\lambda\}$ is a
brsr in $\a$ then we can choose $\{E'_\lambda\}$ also in $\a$.
\end{lem}

\begin{proof}
For simplicity, we assume that $I=[0,\beta]$ ($\alpha=0$).
The general case follows from this by reparametrization. Let
$\lambda_0\in\text{At}(\{E_\lambda\})$, $p_0=p(\lambda_0)$ be the
jump projection at $\lambda_0$ and $\alpha_0=\tau(p_0)$.

It is well known \cite{arvkad,{kad}} that there exists
$\{U_\lambda\}_{\lambda\in [0,\alpha_0]}$ a brsr of $p_0$ in $\M$
such that
\begin{equation}\label{dist uniforme 2}
\tau(U_\lambda)= \frac{\tau(p_0)(\alpha_0-\lambda)}{\alpha_0},\ \
\lambda\in[0,\alpha_0].
\end{equation} Moreover,
 if $\a\subseteq \M$ is a masa and $p_0\in\P(\a)$ then we
can choose $\{U_\lambda\}$ to be in $\a$. Let $$ E'_\lambda=
\left\{
\begin{array}{ll}
    E_\lambda  & \hbox{if }\ 0\leq \lambda<\lambda_0 \\
    E_{\lambda_0}+ U_{\lambda-\lambda_0} & \hbox{if }\ \lambda_0
    \leq \lambda\leq\lambda_0+\alpha_0 \\
    E_{\lambda-\alpha_0} & \hbox{if }\ \lambda_0+\alpha_0<\lambda\leq\alpha_0+\beta. \\
\end{array}
\right. $$ It is easy to see that $\{E'_\lambda\}_{\lambda\in
I'}$, where $I'=[0,\beta+\alpha_0]$, is a brsr. Note that if
$\{E_\lambda\}$ is in a masa $\a\subseteq \M$ then $p_0\in \a$ and
we can choose $\{U_\lambda\}$ in $\a$, so that $\{E'_\lambda\}$ is
also in $\a$.  The increasing, right-continuous function
$f:I\rightarrow I'$ given by
\begin{equation}\label{defi de f}
f(\lambda)=  \left\{
\begin{array}{ll}
\lambda    & \hbox{if }\ 0\leq \lambda<\lambda_0\\
    \lambda+\alpha_0 & \hbox{if }\ \lambda_0\leq
    \lambda\leq\beta+\alpha_0\\
\end{array}
\right. \end{equation} satisfies $E_\lambda=E'_{f(\lambda)}$,
$\lambda\in [0,\beta]$. Moreover
At$(\{E'_\lambda\})=f(\text{At}(\{E_\lambda\})\setminus
\{\lambda_0\})$ and $p(\lambda)=p\,'(f(\lambda))$ for every
$\lambda\in \text{At}(\{E_\lambda\})\setminus \{\lambda_0\}$,
where $p\,'(f(\lambda))$ is the jump projection of
$\{E'_\lambda\}$ at $f(\lambda)\in\text{At}(\{E'_\lambda\})$.
Therefore $$\mathcal J(\{E'_\lambda\})=\sum_{\lambda\in
\text{At}(\{E_\lambda\})\setminus \{\lambda_0\}}
\tau(p(f(\lambda)))=\mathcal J(\{E_\lambda\})-\tau(p_0).$$ The
rest of the properties of $f$ follow directly from (\ref{defi de
f}).
\end{proof}

\medskip

\begin{proof}[Proof of Theorem \ref{prepo}] Let $a\in \M^+$ and
consider the brsr induced by $a$ (see (\ref{la resu de a})). Set
$\beta=\|a\|$, let $I=[0,\beta]$ and let $\{\lambda_n\}_{n\in N}$
be an enumeration of the set At$(\{E_\lambda\})$, where
$N\subseteq \NN$ is an initial segment, and let
$\alpha_n=\tau(p(\lambda_n))>0$. By (\ref{hay a lo sumo
numerables}) we have $\sum_{n\in N}\alpha_n\leq 1$. Let $I_1:=I$,
$\{E^1_\lambda\}:=\{E_\lambda\}$ and let
$(\{E^2_\lambda\}_{\lambda\in I_2},f_1)$ be the strong refinement
obtained from $\{E^1_\lambda\}_{\lambda\in I_1}$ and the atom
$\lambda_1$ as in Lemma \ref{primer refinamiento}. Recall that in
this case $I_2=[0,\beta+\tau(p_1)]$ and set
$g_2:=f_1:I_1\rightarrow I_2$.

 We proceed inductively: assume that for $1\leq t\leq k-1$
  we have
  brsr's
$\{E^t_\lambda\}_{\lambda\in I_t}$, where
$I_t=[0,\beta+\sum_{j=1}^{t-1}\alpha_j]$ and for $1\leq i\leq k-2$
increasing right-continuous functions $f_i:I_{i}\rightarrow
I_{i+1}$ such that $(\{E^{i+1}_\lambda\},f_i)$ strongly refines
$\{E^i_\lambda\}$ and such that $f_i(\lambda)-\lambda\leq
\alpha_i$ for $\lambda\in I_i$.
 Assume further that for $2\leq l\leq k-1$ there exist injective
  functions
$g_l:I\rightarrow I_l$ such that
$$\text{At}(\{E^{\,l}_\lambda\})=g_l(\text{At}(\{E_\lambda\})\setminus
\{\lambda_1,\ldots,\lambda_{l-1}\})$$ and  $$ \mathcal
J(\{E^l_\lambda\})=\mathcal
J(\{E_\lambda\})-\sum_{j=1}^{l-1}\alpha_j.$$ Apply Lemma
\ref{primer refinamiento} to the brsr
$\{E^{k-1}_\lambda\}_{\lambda\in I_{k-1}}$ and the atom
$g_{k-1}(\lambda_{k-1})$. Then we obtain  a brsr
$\{E^{k}_\lambda\}_{\lambda\in I_{k}}$,
$I_{k}=[0,\beta+\sum_{j=1}^{k-1}\alpha_j]$, and an increasing
right-continuous function $f_{k-1}:I_{k-1}\rightarrow I_{k}$ such
that $(\{E^k_\lambda\},f_{k-1})$ is a strong refinement of
$\{E^{k-1}_\lambda\}$; in this case we have
$f_{k-1}(\lambda)-\lambda\leq \alpha_{k-1}$.
If we let
 $g_{k}=f_{k-1}\circ g_{k-1}:I\rightarrow I_{k}$ then $g_{k}$ is injective and
such that
 \begin{eqnarray*}\text{At}(\{E^{k}_\lambda\})&=&f_{k-1}(\text{At}
 (\{E^{k-1}_\lambda\})\setminus\{g_{k-1}(\lambda_{k-1})\})\\ &=&g_{k}( \text{At}
 (\{E_\lambda\})\setminus\{\lambda_1,\ldots,
 \lambda_{k-1}\}).\end{eqnarray*} Moreover,
 $\mathcal J(\{E^{k}_\lambda\})=\mathcal J(\{E^{k-1}_\lambda\})-\alpha_{k-1}=\mathcal J
 (\{E_\lambda\})
 -\sum_{i=1}^{k-1}\alpha_i$.

 We obtain in this way a sequence
$\{E^k_{\lambda}\}_{\lambda\in I_k}$ of brsr's where
$I_k=[0,\beta+\sum_{j=1}^{k-1}\alpha_j]$,
 and increasing right-continuous functions $\{f_k:I_{k}\rightarrow I_{k+1}\}$ for $k\in
\NN$ as in the hypothesis of Lemma \ref{hay limite}.
 Thus, there exists a brsr
$\{E'_\lambda\}_{\lambda\in I'}$ such that for every $k\in \NN$
 $\{E'_\lambda\}$ is a strong refinement of $\{E^k
 _\lambda\}$. In particular, $\{E'_\lambda\}$ is a strong refinement of $\{E
 _\lambda\}=\{E^1_\lambda\}$.
 By Lemma \ref{decrece el salto total}, $\mathcal
J(\{E'_\lambda\})\leq \mathcal J(\{E^k_\lambda\})$ for every $k\in
\NN$ and therefore $\mathcal J(\{E'_\lambda\})=0$, i.e.
$\{E'_\lambda\}$ is diffuse.

Note that if $a\in \cA^+$ for some masa $\a\subseteq \M$ then
$\{E_\lambda\}$ is a brsr in $\a$; by Lemma \ref{primer
refinamiento} we can construct each $\{E_\lambda^k\}$ also in $\a$
and so, by Lemma \ref{hay limite} then $\{E'_\lambda\}$ is in
$\a$. Finally if we let $a'=\int_{I'} \lambda \ dE'_\lambda$ (see
(\ref{int contra resu})) then $a'\in\M^+$ has the desired
properties.
\end{proof}

\section{Modelling of operators and applications}\label{modelado}

\subsection{Modelling of operators}

 We begin with the following elementary lemmas about
functions that we shall need in the sequel.

\begin{lem}\label{funciones decrecientes} Let
$I=[\alpha,\beta],\,J=[\alpha',\beta']\subseteq \RR$ be closed
intervals, $g:J\rightarrow[0,1]$ a decreasing right-continuous
function and let $h:I\rightarrow [0,1]$ be a decreasing continuous
function such that $h(\alpha)\geq g(\alpha')$ and $h(\beta)\leq
g(\beta')$. If we let $\tilde g:J\rightarrow I$ be given by
\begin{equation*}\label{defi de g tilde}
\tilde g(x)=\max\{t\in I:\ g(x)=h(t)\}
\end{equation*}then $\tilde g$ is an increasing right-continuous
function and $g=h\circ \tilde g$.
\end{lem}

\begin{lem}\label{funcion asociada}
Let $I=[\alpha,\beta],\,J=[\alpha',\beta']\subseteq \RR$ and let
$f:J\rightarrow I$ be an increasing right-continuous function such
that $f(\beta')=\beta$.
 If $f^\dagger:I\rightarrow J$ is the function given by
$$f^\dagger(\lambda)=\min\{t\in J:\ \lambda\leq f(t)\}$$ then it
is increasing, left-continuous
  and such that for every $t\in J$
\begin{equation}\label{ecua fun asoc} \{\lambda\in I:\ \lambda>f(t)\}=
\{\lambda\in I:\ f^\dagger(\lambda)> t\}.
\end{equation} If $f$ is strictly increasing then $f^\dagger$ is
continuous. Moreover, if $\tilde J:=[\gamma,\delta]\subseteq J$
and $g:\tilde J\rightarrow I$ is increasing and right-continuous,
$g(\delta)=\beta'$ and
 $f(t)\geq
g(t)$ for every $t\in \tilde J$, then $g^\dagger\geq f^\dagger$.
\end{lem}

\begin{lem}\label{funcion asociada2} Let
$I=[\alpha,\beta],\,J=[\alpha',\beta']\subseteq \RR$ and let
$f:I\rightarrow J$ be an increasing left-continuous function such
that $f(\alpha)=\alpha'$. If $f_\dagger:J\rightarrow I$ is the
function given by $$f_\dagger(\lambda)=\max\{t\in I:\ \lambda\geq
f(t)\}$$ then it is increasing, right-continuous and
 such that for every $t\in I$
\begin{equation}\label{ecua fun asoc2} \{\lambda\in J:\ \lambda <f(t)\}=
\{\lambda\in J:\ f_\dagger(\lambda)< t\}.
\end{equation}
\end{lem}

The following theorem develops the modelling of positive operators
and relates it with the refinement between the spectral
resolutions induced by these operators.

\begin{teo}
\label{primer teo} Let $(\M,\tau)$ be a II$_1$ factor, let $a\in
\M^+$ with continuous distribution and let $I=[0,\|a\|]$. Then
\begin{enumerate}
\item[1.] If $b\in \M^+$, there exists a nonnegative increasing left-continuous
function $h_b$ on $I$ such that if $\tilde b=h_b(a)$ then
$\mu_b=\mu_{\tilde b}$.

\item[2.] The brsr induced by $a$ refines the brsr
induced by $b$ if and only if $\tilde b=b$. Moreover, if the brsr
induced by $a$ strongly refines the brsr induced by $b$ then $h_b$
is continuous.

\item[3.] If $c^+\in \M$ then $c\precsim b$ (resp $c\prec_w b$,
$c\prec b$) if and only if $h_c(a)\leq h_b(a)$ (resp.
$h_c(a)\prec_w h_b(a)$, $h_c(a)\prec h_b(a)$).
\end{enumerate}
\end{teo}

\begin{proof} Let $a\in \M^+$ with continuous
distribution, let $I=[0,\|a\|]$ and let $h:I\rightarrow [0,1]$ be
the decreasing continuous function defined by
$h(t)=\tau(P^a(t,\infty))$. Note that $h(\|a\|)=0$ and, since $a$
has continuous distribution, $h(0)=1$.

 Let $b\in \M^+$, $J=[0,\|b\|]$ and let $g:J\rightarrow [0,1]$ be the
decreasing right-continuous function defined by
$g(s)=\tau(P^b(s,\infty))$. By Lemma \ref{funciones decrecientes},
there exists an increasing right-continuous function $\tilde
g:J\rightarrow I$, such that $g=h\circ \tilde g$, i.e.
\begin{equation}\label{igual de traza}
\tau(P^b(s,\infty))=\tau(P^{a}(\tilde g(s),\infty)),\ \ \ s\in J.
\end{equation}
 By Lemma \ref{funcion asociada} there exists an
increasing (and therefore uniformly bounded measurable)
left-continuous function $h_b:=\tilde g^\dagger:I\rightarrow J$
such that
\begin{equation}\label{ecua para b} \{\lambda\in I:\
h_b(\lambda)> s\}=\{\lambda\in I:\ \lambda>\tilde g (s)\}, \ \
s\in J.
\end{equation}
Let $\tilde b= h_b(a)$ and note that $\tau(P^{\tilde
b}(s,\infty))=\tau(P^{b}(s,\infty))$, which follows from
(\ref{igual de traza}) and (\ref{ecua para b}). Therefore, $b$ and
$\tilde b$ have the same singular values.

To prove {\it 2.} assume that the brsr induced by $b\in \M^+$
is refined by the brsr induced by $a$. Let $\tilde b=h_b(a)$ and
note that $P^{\tilde b}(s,\infty)= P^{a}(\tilde
g(s),\infty)$ and by hypothesis $P^b(s,\infty)=P^{a}(f(s),\infty)$
for some increasing right-con\-tinuous function $f:J\rightarrow
I$. Then $P^{\tilde b}(s,\infty)\leq P^{b}(s,\infty)$ or
$P^{b}(s,\infty)\leq P^{\tilde b}(s,\infty)$ and by (\ref{ecua
para b}) we have $\tau(P^b(s,\infty))=\tau(P^{\tilde
b}(s,\infty))$ so $P^b(s,\infty)=P^{\tilde b}(s,\infty)$, $s\in
J.$ Therefore $b=\tilde b$. On the other hand, if $b=j(a)$ for any
increasing left-continuous function $j:I\rightarrow J$, then by
Lemma \ref{funcion asociada2} there exists an increasing
right-continuous function $f:=j_\dagger:J\rightarrow I$ such that
\begin{eqnarray*}
 P^{b}(\lambda,\infty)&=&P^{a}(\{
t\in I:\ \lambda< j(t)\})\\ &=& P^{a}(\{ t\in I:\ f(\lambda)<
t\})=P^{a}(f(\lambda),\infty),\end{eqnarray*} so the brsr induced
by $a$ refines the brsr induced by $b$. Finally assume that the
brsr induced by $a$ strongly refines the brsr induced by $b$.
Then, by (d) in Definition \ref{defi de refinamiento} $f$ is
strictly increasing and therefore, by Lemma \ref{funcion asociada}
$h_b=f^\dagger $ is continuous.

To prove {\it 3.} assume that $c\in \M^+$ is such that
$\tau(P^c(s,\infty))\leq \tau(P^b(s,\infty))$ for all $s\geq 0$
and therefore $\|c\|\leq\|b\|$. As before, let
$k:[0,\|c\|]\rightarrow [0,1]$ be the function given by
$k(s)=\tau(P^c(s,\infty))$, $\tilde k$ obtained from $k$ as in
Lemma \ref{funciones decrecientes}, and $h_c=\tilde k^\dagger$
obtained from $\tilde k$ as in Lemma \ref{funcion asociada}. Then,
  $\tilde g(t)\leq \tilde k(t)$ for every
$t\in [0,\|c\|]$ and, by Lemma \ref{funcion asociada}, we conclude
that $h_c=\tilde k^\dagger\leq \tilde g^\dagger=h_b$. From this it
follows that $\tilde c\leq \tilde b$, where $\tilde b=h_b(a), \,
\tilde c=h_c(a)$.
 The rest of the statement is a
consequence of the fact that $\mu_b=\mu_{\tilde b}$ and
$\mu_c=\mu_{\tilde c}$.
\end{proof}

 We say that $c\in \M^+$ is a \emph{model} of $b\in \M^+$
with respect to $a\in \M^+$, if there exists a nonnegative,
left-continuous and increasing function $h$ such that $c=h(a)$ and
$\mu_c=\mu_b$. Thus, with the notations of the proof of Theorem
\ref{primer teo}, we see that $\tilde b\in \M^+$ is a model of
$b\in \M^+$ with respect to $a$. As an immediate consequence of
{\it 2.} in Proposition \ref{rems}, we conclude that the model
$\tilde b$ is approximately unitarily equivalent to $b$ in $\M$.

\begin{rem} In \cite{kad} Kadison solved
the following problem in a II$_1$ factor $(\M,\tau)$: given a masa
$\cA\subseteq \M$, $a\in \cA_{sa}$ and $t\in [0,1]$ find a
projection $p\in\cA$ and $\lambda\in \RR$ such that $\tau(p)=t$,
$ap\geq \lambda p$ and $a(I-p)\leq \lambda(I-p)$. Note that
Theorems \ref{prepo} and \ref{primer teo} give an alternative
proof of this statement in the case $a\in \cA^+$. Indeed, let
$a'\in\cA^+$ be as in Theorem \ref{prepo} and $h_a$ be as in
Theorem \ref{primer teo}. Then, if we let
$p=P^{a'}(\alpha,\infty)$ with $\tau(p)=t$ (such $\alpha$ always
exists since $a'$ has continuous distribution) and
$\lambda=h_a(\alpha)$ then $p$ and $\lambda$ have the desired
properties, since $h_a$ is an increasing function. \qed
\end{rem}

 As a final comment let us note that a variation of the
proof of Theorem \ref{primer teo} implies that if $a\in \M^+$ has
continuous distribution and if $\nu$ is any regular Borel
probability measure of compact support in the real line then,
there exists $h:[0,\|a\|]\rightarrow \RR$ such that
$\nu(\Delta)=\tau(P^{h(a)}(\Delta))$.
 Indeed, we just have to replace the function
 $\tau(P^b(\lambda,\infty))$ by $\nu((\lambda,\infty))$ in the
 proof of {\it 2.}
In particular, if $\cA\subseteq \M$ is a masa and we consider
$a\in \cA^+$ then this argument gives a different proof of
Proposition 5.2 in \cite{arvkad}.

\subsection{Some applications of the modelling technique}

The following application of Theorem \ref{primer teo} provides new
characterizations of spectral preorder and sub-majorization
between positive operators in II$_1$ factors. Note that these
re-formulations have an inequality-type form.

\begin{teo}\label{segundo teo} Let $(\M,\tau)$ be a II$_1$
factor and let $a,\,b\in \M^+$.
Then
\begin{enumerate}
\item[1.] $b$ spectrally
dominates $a$ if and only if
\begin{equation}\label{primer carac}
\peso{there exists} c \in \cluni{b} \peso{with} a\leq c
\end{equation}
or, equivalently, if
\begin{equation}\label{segunda carac}
 \peso{there exists} d \in \cluni{a} \peso{with} d\leq b.
\end{equation}
Moreover, we can assume that $a$ and $c$ commute and that $b$ and
$d$ commute.
\item[2.] $b$ sub-majorizes $a$  if and only if there exists $c\in
\M^+$ such that \begin{equation}\label{carac de la submayo} a\leq
c\prec b. \end{equation} Moreover, we can assume that $a$ and $c$
commute.
\end{enumerate}
\end{teo}
\begin{proof}
Recall that for positive operators $a,\,b\in \M^+$, $a\leq b$
implies $a\precsim b$. Thus, the existence of a sequence of
unitary operators satisfying (\ref{primer carac}) or (\ref{segunda
carac}) implies spectral domination. Analogously, the existence of
an operator satisfying (\ref{carac de la submayo})  implies
sub-majorization. Next show that the reverse implications are also
true.

To prove the first part of {\it 1.} let $a,\,b\in\M^+$ such that
$a\prec b$. By Theorem \ref{prepo} there exists $a'\in\M^+$ with
continuous distribution such that the brsr induced by $a'$
(strongly) refines the brsr induced by $a$. By Theorem \ref{primer
teo} there exists an increasing left-continuous function $h_b$
such that, if $\tilde b=h_b(a')$, $\mu_b=\mu_{\tilde b}$. By {\it
2.} in Proposition \ref{rems}, this implies that $\tilde b \in
\cluni{b}$. Since by hypothesis $\mu_a\leq \mu_b$, by {\it 2.} and
{\it 3.} in Theorem \ref{primer teo} we have $\tilde b=h_b(a')\geq
h_a(a')=a$. Thus, we obtain (\ref{primer carac}) with
$c=\tilde b$. The proof of the second part follows a similar path,
considering the model of $a$ with respect to a refinement of $b$.

To prove {\it 2.}, let $a$ and $a'$ be as in the first part of the
proof. Let $b\in \M^+$ be such that $a\prec_w b$ and let $\nu$
denote the regular Borel probability measure on $I'=[0,\|a'\|]$
given by $\nu(\Delta)=\tau(P^{a'}(\Delta))$. Then, if $h_a,\,h_b$
are as in Theorem \ref{primer teo} we have (see Remark
\ref{fundamental sobre funciones}) that $h_a\prec_w h_b$ in
$L^\infty(\nu)$.
 Therefore, by Proposition \ref{pro HiaiNak} there
exists $h\in L^\infty(\nu)$ such that $h_a\leq h\prec h_b$. If we
let $c=h(a')$ then $a\leq c\prec b$ by construction,
since $a=h_a(a')$.
\end{proof}

The first part of {\it 1.} in Theorem \ref{segundo teo} gives a
partial affirmative solution to the following problem posed in
 \cite{{Doug1},Doug}: given a $(\M,\tau)$ a II$_1$ factor and
 $a,\,b\in\M^+$ such that $a\precsim b$, is there any automorphism
 of $\M$, $\Theta$, such that $\Theta(b)\geq a$? Our
 considerations above lead to a sequence of $\tau$-preserving
 automorphisms (Ad$_{u_n}$)$_{n\in \NN}$, where $u_n\in \U_\M$, such
 that \emph{in the limit} the above statement is true.

\begin{cor}\label{mas caracterizaciones}
Let $a,\,b\in \M^+$. Then the following statements are equivalent:
\begin{enumerate}
\item[1.] $b$ spectrally dominates $a$.
\item[2.] There exists a brsr
$\{E_\lambda\}_{\lambda\in I}$, where $I=[0,\|a\|]$ such that
$\tau(E_\lambda)=\tau(P^a(\lambda,\infty))$ for every $\lambda\in
I$  and
\begin{equation}\label{con la reso} \lambda E_\lambda\leq E_\lambda \,b\,
E_\lambda, \ \ \ \forall \lambda\geq 0.\end{equation}
\end{enumerate}
\end{cor}

\begin{proof}
Assume {\it 1.} and note that, by Theorem \ref{segundo teo} there
exists a sequence $(v_n)_n\subseteq \U_\M$ such that
$\lim_{n\rightarrow\infty}\|d-v_n^*av_n\|=0$ and $d\leq b$ for
some $d\in \M^+$. Then $\tau(p(a))=\tau(p(d))$ for every
polynomial $p\in \CC[x]$ and, using monotone convergence,
 we have
$\tau(P^a(\lambda,\infty))=\tau(P^d(\lambda,\infty))$,
$\lambda\geq 0$.
 Moreover, $$\lambda P^d(\lambda,\infty)\leq
P^d(\lambda,\infty)\,d\leq
P^d(\lambda,\infty)\,b\,P^d(\lambda,\infty).$$ Then, if we set
$E_\lambda=P^d(\lambda,\infty)$,
$\{E_\lambda\}_{\lambda\in[0,\,\|a\|]}$ is the desired brsr.
Conversely, assume that there exists a brsr
$\{E_\lambda\}_{\lambda\in [0,\|a\|]}$ as in item 2. Given
$\epsilon
>0$, let $b_\epsilon=b+\epsilon I$ and note that
$\lambda\,E_\lambda<E_\lambda b_\epsilon E_\lambda$, so
$P^{E_\lambda \,b_\epsilon\,E_\lambda}(\lambda,\infty)=E_\lambda$.
In \cite{fack} Fack proved the following inter\-lacing-like
inequality: for every orthogonal projection $p\in \M$,
$p\,b\,p\precsim b$. Then we have $$\tau(P^a(\lambda,\infty))=
\tau(E_\lambda)=\tau(P^{E_\lambda
\,b_\epsilon\,E_\lambda}(\lambda,\infty))\leq
\tau(P^{b_\epsilon}(\lambda,\infty)).$$ The inequality above shows
that $\mu_a\leq \mu_{b_\epsilon}$ for every $\epsilon>0$. The
corollary is now a consequence of the fact that
$\lim_{\epsilon\rightarrow 0^+}\mu_{b_\epsilon}(t)=\mu_b(t)$ for
every $t\geq 0$.
\end{proof}

We end with some applications of our previous results. These are
mostly re-statements of some inequalities with respect to spectral
preorder and sub-majorization obtained in
\cite{JD,{arvkad},silva,Doug}, using Theorem \ref{segundo teo}.

\begin{cor} \label{algunas desi de young y jensen}Let $(\M,\tau)$ be a II$_1$ factor.
\begin{enumerate}
\item[1.] (Young-type inequalities) Let $x,\,y\in \M$ and
let $p,\,q$ be conjugated indices. Then there exist sequences
$(u_n)_{n\in \NN},(v_n)_{n\in \NN}\subseteq \U_\M$ such that
 $$|xy^*|\leq\lim_{n\rightarrow \infty} u_n^*
(p^{-1}|x|^p+q^{-1}|y|^q)u_n$$ and $$ \lim_{n\rightarrow \infty}
v_n^*|xy^*| v_n\leq p^{-1}|x|^p+q^{-1}|y|^q$$
\item[2.] (Jensen-type inequalities) Let $\cA$ be a unital
$C^*$-algebra, $\Phi:\cA\rightarrow\M$ be a positive unital map,
$a\in \cA^+$ and $f:\sigma(a)\rightarrow \RR$ be a convex
function.
\begin{enumerate}
\item[(a)]  If $f$ is increasing, 
there exist sequences $(u_n)_{n\in \NN},\,(v_n)_{n\in
\NN}\subseteq  \U_\M$ with $$f(\Phi(a))\leq
\lim_{n\rightarrow \infty} u_n^*\Phi(f(a))u_n$$ and
$$\lim_{n\rightarrow \infty} v_n^* f(\Phi(a))v_n\leq \Phi(f(a)).$$
\item[(b)] If $f$ is an arbitrary convex function, there exists
$c\in \M^+$ such that $$f(\Phi(a))\leq c\prec \Phi(f(a)).$$
Moreover, we can choose $c$ so that it commutes with $f(\Phi(a))$.
\end{enumerate}

\end{enumerate}
\end{cor}
\begin{proof}
In \cite{Doug}, Farenick and Manjegani proved that if
$p,\,q,\,x,\,y$ are as above, then $|xy^*|\precsim
p^{-1}|x|^p+q^{-1}|y|^q$. On the other hand, in \cite{JD} it was
shown that if $\Phi,\,f,\,a$ are as above then,
$f(\Phi(a))\precsim\phi(f(a))$ if $f$ is increasing and in
general, $f(\Phi(a))\prec_w \Phi(f(a))$ for an arbitrary convex
function $f$.
 The corollary
follows from these facts and Theorem \ref{segundo teo}.
\end{proof}

The proofs of Theorem \ref{segundo teo} and Corollary \ref{mas
caracterizaciones} show a possible interplay between Theorems
\ref{prepo} and \ref{primer teo} to get an interesting tool to
deal with problems regarding spectral relations. As far as we
know, the conclusions of Corollary \ref{algunas desi de young y
jensen} are not possible using the previous literature.

Some of our results extend to certain classes of
(unbounded) measurable operators affiliated with $\M$. Also, note
that there is still the problem of finding characterizations of
spectral order and sub-majorization similar to those in Theorem
\ref{segundo teo}, for general semifinite factors; these
characterizations may depend on generalizations  of both Theorems
\ref{prepo} and \ref{primer teo}. We shall investigate these
matters elsewhere.

\medskip

\noindent{\bf Acknowledgements.} I would like to thank Professors
Demetrio Stojanoff, Doug Farenick and Martin Argerami for fruitful
discussions on the material contained in this note. I would also
like to thank the referee for several suggestions that improved
the paper.

\bigskip

{\noindent  \small Pedro G. Massey \\ Departamento de
Matem\'atica, FCE-Universidad Nacional de La Plata and \\ Instituto Argentino de Matemática - CONICET, Argentina\\
Mail address: CC 172, (CP 1900) La Plata, Prov. de Buenos Aires,
Argentina
\\ E-mail address: massey@mate.unlp.edu.ar}


\begin{thebibliography}{99}

\bibitem{JD} J. Antezana, P. Massey, and D. Stojanoff,
{\it Jensen's Inequality and Majorization}, pre\-print ({\tt
http://xxx.lanl.gov/abs/math.FA/0411442}).

\bibitem{argmas} M. Argerami and P. Massey, {\it The local form of doubly stochastic
maps and joint majorization in II$_1$ factors}, preprint.

\bibitem{argmas2} M. Argerami and P. Massey, {\it On a Schur-Horn
type theorem in II$_1$ factors}, preprint ({\tt
http://arxiv.org/abs/math.OA/0604120}).

\bibitem{arvkad}W. Arveson and R. Kadison, {\it Diagonals of
self-adjoint operators}, preprint ({\tt
http://xxx.lanl.gov/abs/math.OA/0508482}).

\bibitem{silva} J.S. Aujla and F.C. Silva, {\it Weak
majorization inequalities and convex functions}, Linear Algebra
Appl. 369 (2003), 217-233.

\bibitem{Doug1} D.R. Farenick, private communication.

\bibitem{Doug} D.R. Farenick and S.M. Manjegani, {\it Young's Inequality in Operator Algebras},
to appear in J. of the Ramanujan Math. Soc. ({\tt
http://xxx.lanl.gov/abs/math.OA/0303318}).

\bibitem{fack} T. Fack, {\it Sur la notion de valeur
caractéristique}, J. Operator Theory (1982), 307-333.

\bibitem {Hiai0} F. Hiai, {\it Majorization and Stochastic maps
in von Neumann algebras}, J. Math. Anal. Appl. {\bf 127} (1987),
no. 1, 18--48.

\bibitem {Hiai} F. Hiai, {\it Spectral majorization between normal
operators in von Neumann algebras}, Operator algebras and operator
theory (Craiova, 1989), 78--115, Pitman Res. Notes Math. Ser.,
271, Longman Sci. Tech., Harlow, 1992.

\bibitem{HiaiN0} F. Hiai, Y. Nakamura, {\it Distance between
unitary orbits in von Neumann algebras}, Pacific J. of Math, vol.
138 (1989) 259-294.

\bibitem{HiaiN} F. Hiai, Y. Nakamura, {\it Closed Convex Hulls of Unitary Orbits in von Neumann Algebras},
Trans. Amer. Math. Soc. {\bf 323} (1991), 1-38.

\bibitem{kad1} R.V. Kadison, {\it The Pythagorean theorem I: the finite
case}. Proc. N.A.S. (USA),99(7):4178-4184, 2002.

\bibitem{kad2}R.V. Kadison, {\it The Pythagorean theorem II: the infinite
case}. Proc. N.A.S. (USA),99(8):5217-5222, 2002.

\bibitem{kad} R. V. Kadison, {\it Non-commutative
conditional expectations and their applicactions}, volume 365 of
Contemporary Mathematics, pages 143-179. Amer. Math Soc. 2004.
Operator algebras, quantization and non-commutative geometry.

\bibitem{kam} E. Kamei,{\it Majorization in finite factors},
Math. Japonica 28, No. 4 (1983), 495-499.

\bibitem{she}D. Sherman, {\it Unitary orbits
of normal operators in von Neumann algebras}, to appear in Crelle's Journal.

\end{thebibliography}
\end{document}